\documentclass{amsart}
\usepackage{amssymb, amsmath, amsfonts, amscd}

\input xypic

\theoremstyle{plain}
\newtheorem{theorem}{Theorem}[section]
\newtheorem{corollary}[theorem]{Corollary}
\newtheorem{lemma}[theorem]{Lemma}
\newtheorem{proposition}[theorem]{Proposition}

\newtheorem{example}[theorem]{Example}

\theoremstyle{definition}
\newtheorem{definition}[theorem]{Definition}

\theoremstyle{remark}

\numberwithin{equation}{theorem}

\newcommand{\End}{\operatorname{End}}

\newcommand{\F}{\mathcal{F}}

\newcommand{\I}{\mathcal{I}}

\renewcommand{\L}{\mathcal{L}}
\newcommand{\E}{\mathcal{E}}
\newcommand{\V}{\mathbb{V}_{\E}}
\newcommand{\VL}{\mathbb{V}_{\L}}

\newcommand{\Hom}{\operatorname{Hom}} 
\renewcommand{\O}{\mathcal{O} }

\renewcommand{\Pr}{\mathcal{J} }

\newcommand{\pa}{\partial}

\newcommand{\KS}{\operatorname{KS}}
\newcommand{\dlog}{\operatorname{dlog}}
\newcommand{\AT}{\operatorname{AT}}
\newcommand{\LR}{\operatorname{LR}}

\renewcommand{\H}{\operatorname{H} }
\newcommand{\HH}{\operatorname{HH}}
\newcommand{\Pic}{\operatorname{Pic} }

\newcommand{\exan}{\operatorname{exan}}
\newcommand{\Exan}{\operatorname{Exan}}

\newcommand{\Ext}{\operatorname{Ext} }

\newcommand{\Der}{\operatorname{Der} }

\newcommand{\one}{\mathbf{1}}
\newcommand{\spec}{\operatorname{Spec}}
\newcommand{\id}{\operatorname{id}}

\renewcommand{\one}{\mathbf{1}}

\begin{document}

\title{On jets, extensions  and characteristic classes I}

\author{Helge Maakestad}
\address{Oslo, Norway}

\email{\text{h\_maakestad@hotmail.com} }
\keywords{Atiyah sequence, jet bundle, Lie-Rinehart algebra, Chern
  class, Atiyah class, Kodaira-Spencer class, square zero extension, lifting}
\thanks{Partially supported by a scholarship from NAV, www.nav.no}

\subjclass{14F10, 14F40}

\date{Spring 2009}

\begin{abstract}  In this paper we give general definitions
of non-commutative jets in the local and global situation
using square zero extensions and derivations. We study the functors
$\Exan_k(A,I)$ where $A$ is any $k$-algebra and $I$ is any left and
right $A$-module and use this to relate affine non-commutative jets to
liftings of modules. We also study the Kodaira-Spencer class $\KS(\L)$
and relate it to the Atiyah class.
\end{abstract}

\maketitle

\tableofcontents

\section{Introduction} 

In this paper we give general definitions of non-commutative
jets in the local and global situation using square zero extensions
and derivations.
We study the functors $\Exan_k(A,I)$ where $A$ is any $k$-algebra and
$I$ is any left and
right $A$-module and use this to relate affine non-commutative jets to
liftings of modules.
In the final section of the paper we define and prove basic properties
of the Kodaira-Spencer class $\KS(\L)$ and relate it to the Atiyah
class.

\section{Jets, liftings and small extensions}

We give an elementary discussion of structural properties of square
zero extensions of arbitrary associative unital $k$-algebras. We
introduce for any $k$-algebra $A$ and any left and right $A$-module
$I$ the set $\Exan_k(A,I)$ of isomorphism classes of square zero
extensions of
$A$ by $I$ and show it is a left and right module over the center
$C(A)$ of $A$. This structure generalize the structure as left
$C(A)$-module introduced in \cite{grothendieck}. We also give an
explicit construction of $\Exan_k(A,I)$ in terms of cocycles. Finally
we give a direct construction of non-commutative jets and generalized
Atiyah sequences using derivations and square zero extensions.

Let in the following $k$ be a fixed base field and let
\[ 0\rightarrow I \rightarrow B \rightarrow A \rightarrow 0 \]

be an exact sequence of associative unital $k$-algebras with
$i(I)^2=0$. Let $i:I\rightarrow B$ and $p:B\rightarrow A$ denote the
morphisms.
Assume $s$ is a
map of $k$-vector spaces with the following properties:
\[ s(1)=1\]
and
\[ p\circ s=\id.\]
Such a section always exist since $B$ and $A$ are vector spaces over
the field $k$. Note: $s$ gives the ideal $I$ a left and right
$A$-action.

\begin{lemma} There is an isomorphism
\[ B\cong I\oplus A\]
of $k$-vector spaces.
\end{lemma}
\begin{proof} Define the following maps of vector spaces:
$ \phi: B\rightarrow I\oplus A$ by $\phi(x)=(x-sp(x),p(x))$
and $\psi:I\oplus A\rightarrow B$ by $\psi(u,x)=u+s(x).$
It follows $\psi\circ \phi=\id$ and $\phi\circ \psi=\id$
and the claim of the Proposition follows.
\end{proof}

Define the following element:
\[ \tilde{C}:A\times A\rightarrow I\]
by
\[ \tilde{C}(x\times y)=s(x)s(y)-s(xy).\]

It follows that $\tilde{C}=0$ if and only if $s$ is a ring
homomorphism.

\begin{lemma} The map $\tilde{C}$ gives rise to an element $C\in
  \Hom_k(A\otimes_k A,I)$.
\end{lemma}
\begin{proof}
We easily see that $\tilde{C}(x+y,z)=\tilde{C}(x,z)+\tilde{C}(y,x)$
and
$\tilde{C}(x,y+z)=\tilde{C}(x,y)+\tilde{C}(x,z)$ for all $x,y,z\in
A$. Moreover for any $a\in
k$ it follows
\[ \tilde{C}(ax,y)=\tilde{C}(x,ay)=a\tilde{C}(x,y).\]
Hence we get a well defined element $C\in \Hom_k(A\otimes_k A,I)$ as
claimed.
\end{proof}

Define the following product on $I\oplus A$:
\begin{align}
\label{prod} (u,x)\times (v,y)=&(uy+xv+C(x,y),xy).
\end{align}
We let $I\oplus^C A$ denote the abelian group $I\oplus A$ with product
defined by \ref{prod}.

\begin{proposition} \label{main} The natural isomorphism
\[B\cong I\oplus A\]
of vector spaces is a unital ring isomorphism if and only if the
following holds:
\[ xC(y,z)-C(xy,z)+C(x,yz)-C(x,y)z=0\]
for all $x,y,z\in A$.
\end{proposition}
\begin{proof} We have defined two isomorphisms of vector spaces
  $\phi,\psi$:
\[ \phi(x)=(x-sp(x),p(x))\]
and
\[ \psi(u,x)=u+s(x).\]
We define a product on the direct sum $I\oplus A$ using $\phi$ and
$\psi$:
\[ (u,x)\times
(v,y)=\phi(\psi(u,x)\psi(v,y))=\phi((u+s(x))(v+s(y)))=\]
\[ \phi(uv+us(y)+s(x)v+s(x)s(y))=\]
\[(us(y)+s(x)v+s(x)s(y)-s(xy),xy)=(uy+xv+C(x,y),xy).\]
Here we define
\[ uy=us(y)\]
and
\[xv=s(x)v.\]

One checks that
\[ \phi(1)=(1-sp(1),1)=(0,1)=\mathbf{1}\]
and
\[ \one (u,x)=(u,x)\one=(u,x)\]
for all $(u,x)\in I\oplus A$. It follows the morphism $\phi$ is
unital.
Since $C(x+y,z)=C(x,z)+C(y,z)$ and $C(x,y+z)=C(x,y)+C(x,z)$ the
following holds:
\[ (u,x)((v,y)+(w,z))=(u,x)(v,y)+(u,x)(w,z) \]
and
\[ ((v,y)+(w,z))(u,x)=(v,y)(u,x)+(w,z)(u,x).\]

Hence the multiplication is distributive over addition.
Hence for an arbitrary section $s$ of $p$ of vector spaces mapping the
identity to the identity it follows the multiplication defined above
always has a left and right unit and is distributive.
We check when the multiplication is associative.
\[ ((u,x)(v,y))(w,z)=(uyz+xvz+xyw+C(x,y)z+C(xy,z),xyz).\]
Also
\[(u,x)((v,y)(w,z))=(uyz+xvz+xyw+xC(y,z)+C(x,yz),xyz).\]
It follows the multiplication is associative if and only if the
following equation holds for the element $C$:
\[ xC(y,z)-C(xy,z)+C(x,yz)-C(x,y)z=0\]
for all $x,y,z\in A$. The claim follows.
\end{proof}

Let
\begin{align}
\label{eq1} xC(y,z)-C(xy,z)+C(x,yz)-C(x,y)z=0.
\end{align}
be the \emph{cocycle condition}.

\begin{definition} Let $\exan_k(A,I)$ be the set of elements $C\in
  \Hom_k(A\otimes_k A,I)$ satisfying the cocycle condition \ref{eq1}.
\end{definition}

\begin{proposition}\label{cocycle} Equation \ref{eq1} holds for all
  $x,y,z\in A$:
\end{proposition}
\begin{proof} We get:
\[ xC(y,z)=s(x)s(y)s(z)-s(x)s(yz).\]
\[ C(xy,z)=s(xy)s(z)-s(xyz).\]
\[ C(x,yz)=s(x)s(yz)-s(xyz),\]
and
\[C(x,y)z=s(x)s(y)s(z)-s(xy)s(z).\]
We get
\[ xC(y,z)-C(xy,z)+C(x,yz)-C(x,y)z= s(x)s(y)s(z)-s(x)s(yz)-\]
\[s(xy)s(z)+s(xyz)+s(x)s(yz)-s(xyz)-s(x)s(y)s(z)+s(xy)s(z)=0\]
and the claim follows.
\end{proof}

\begin{corollary} The  morphism $\phi:B\rightarrow I\oplus^C A$ is an
  isomorphism of unital associative $k$-algebras.
\end{corollary}
\begin{proof} This follows from Proposition \ref{cocycle} and
  Proposition
  \ref{main}.
\end{proof}

Hence there is always a commutative diagram of exact sequences

\[
\diagram 0 \rto & I \rto \dto^= & B \rto \dto^\cong & A \rto \dto^= &
0 \\
         0 \rto & I \rto^i & I\oplus^C A \rto^p & A \rto & 0
\enddiagram
\]
where the middle vertical morphism is an isomorphism associative
unital $k$-algebras.

Define the following left and right $A$-action on the ideal $I$:
\[ xu=s(x)u, ux=us(x) \]
where $s$ is the section of $p$ and $x\in A$, $u\in I$. Recall
$I^2=0$.

\begin{proposition} The actions defined above give the ideal $I$ a
left and right $A$-module structure. The structure is independent of
choice of section $s$.
\end{proposition}
\begin{proof} One checks that for any $x,y\in A$ and $u,v\in I$ the
  following holds:
\[(x+y)u=xu+yu, x(u+v)=xu+xv , 1u=1.\]
Also
\[ (xy)u-x(yu)=s(xy)u-s(x)s(y)u=(s(xy)-s(x)s(y))u=0\]
since $I^2=0$. It follows $(xy)u=x(yu)$ hence $I$ is a left
$A$-module.
A similar argument prove $I$ is a right $A$-module.
Assume $t$ is another section of $p$. It follows
\[ s(x)u-t(x)u=(s(x)-t(x))u=0\]
since $I^2=0$. It follows $s(x)u=t(x)u$. Similarly $us(x)=ut(x)$ hence
$s$ and $t$ induce the same structure of $A$-module on $I$ and the
Proposition is proved.
\end{proof}

We have proved the following Theorem:
Let $A$ be any associative unital $k$-algebra and let $I$ be a left
and right $A$-module. Let $C:A\otimes_kA\rightarrow I$ be a morphism
satisfying the cocycle condition \ref{eq1}.

\begin{theorem} \label{classification} The exact sequence
\[0\rightarrow I\rightarrow I\oplus^C A\rightarrow A\rightarrow 0\]
is a square zero extension of $A$ with the module $I$. Moreover any
square zero extension of $A$ with $I$ arise this way for some
morphism $C\in \Hom_k(A\otimes_k A,I)$ satisfying Equation \ref{eq1}.
\end{theorem}
\begin{proof} The proof follows from the discussion above.
\end{proof}
Let
\[ 0\rightarrow I\rightarrow  E  \rightarrow A \rightarrow 0 \]
with $i:I\rightarrow E$ and $p:E\rightarrow A$
and
\[ 0\rightarrow J \rightarrow  F \rightarrow  B \rightarrow 0 \]
with $j:J\rightarrow F$ and $q:F\rightarrow B$
be square zero extensions of associative $k$-algebras $A,B$ with left
and right modules $I,J$. This
means the sequences are exact and the following holds:
$i(I)^2=j(J)^2=0$.
A triple $(w,u,v)$ of maps of $k$-vector spaces giving rise to a
commutative diagram of
exact sequences
\[
\diagram 0 \rto & I \rto^i \dto^w & E \rto^p \dto^u & A \rto \dto^v &
0 \\
         0 \rto &  J \rto^j        & F \rto^q        & B \rto        &
         0
\enddiagram
\]
is a morphism of extensions if $u$ and $v$ are maps of $k$-algebras
and $w$ is a map of left and right modules. This means
\[ w(x+y)=w(x)+w(y), w(ax)=v(a)w(x), w(xa)=w(x)v(a) \]
for all $x,y\in I$ and $a\in A$.

We say two square zero extensions
\[ 0\rightarrow I \rightarrow  E \rightarrow A \rightarrow 0 \]
and
\[ 0\rightarrow I \rightarrow  F \rightarrow  A \rightarrow 0 \]
are \emph{equivalent} if there is an isomorphism $\phi:E\rightarrow F$
of $k$-algebras making all diagrams commute.

\begin{definition} Let $\Exan_k(A,I)$ denote the set of all
  isomorphism classes of square zero extensions of $A$ by $I$.
\end{definition}

\begin{theorem} \label{module} Let $C(A)$ be the center of $A$.
The set $\exan_k(A,I)$ is a left and right module over
$C(A)$. Moreover there is a bijection
\[ \Exan_k(A,I)\cong \exan_k(A,I) \]
of sets.
\end{theorem}
\begin{proof} We first prove $\exan_k(A,I)$ is a left and right
  $C(A)$-module. Let $C, D\in \exan_k(A,I)$. This means $C, D\in
  \Hom_k(A\otimes_kA ,I)$ are elements satisfying the cocycle
  condition \ref{eq1}. let $a,b\in C(A)\subseteq A$ be elements.
Define $aC,Ca$ as follows:
\[ (aC)(x,y)=aC(x,y)\]
and
\[ (Ca)(x,y)=C(x,y)a.\]
We see
\[ x(aC)(y,x)-(aC)(xy,z)+(aC)(x,yz)-(aC)(x,y)z=\]
\[ a(xC(y,z)-C(xy,z)+C(x,yz)-C(x,y)z)=a(0)=0\]
hence $aC\in \exan_k(A,I)$. Similarly one proves $Ca\in \exan_k(A,I)$
hence we have defined a left and right action of $C(A)$ on the set
$\exan_k(A,I)$.
Given $C,D\in \exan_k(A,I)$ define
\[ (C+D)(x,y)=C(x,y)+D(x,y).\]
One checks that $C+D\in \exan_k(A,I)$ hence $\exan_k(A,I)$ has an
addition operation. One checks the following hold:

\[ a(C+D)=aC+aD, (C+D)a=Ca+Da ,\]
\[(a+b)C=aC+bC, C(a+b)=Ca+Cb ,\]
\[ a(bC)=(ab)C, C(ab)=(Ca)b, 1C=C1=C,\]
hence the set $\exan_k(A,I)$ is a left and right $C(A)$-module.
Define the following map: Let $[B]=[I\oplus^C A]\in \Exan_k(A,I)$ be
an equivalence class of a square zero extension.
Define
\[ \phi:\Exan_k(A,I)\rightarrow \exan_k(A,I) \]
by
\[ \phi[B]=\phi[I\oplus^C A]=C.\]
We prove this gives a well defined map of sets: Assume $[I\oplus^C A]$
and $[I\oplus^D A]$ are two elements in $\Exan_k(A,I)$.
Note: We use brackets to denote isomorphism classes of extensions. The
two
extensions are equivalent if and only if there is an isomorphism
\[ f:I\oplus^C A\rightarrow I\oplus^D A \]
of $k$-algebras such that all diagrams are commutative. This means
\[ f(u,x)=(u,x) \]
for all $(u,x)\in I\oplus^C A$. We get
\[ f((u,x)(v,y))=f(u,x)f(v,y).\]
This gives the equality
\[ (uy+xv+C(x,y),xy)=(uy+xv+D(x,y),xy) \]
for all $(u,x),(v,y)\in I\oplus^C A$. Hence
$\phi[I\oplus^CA]=C=D=\phi[I\oplus^D A]$ and the map $\phi$
is well defined. It is clearly an injective map.
It is surjective by Theorem \ref{classification} and the claim of the
Theorem follows.
\end{proof}

Theorem \ref{module} shows there is a structure of left and right
$C(A)$-module
on the set of equivalence classes of extensions $\Exan_k(A,I)$. The
structure as left $C(A)$-module agrees with the one defined in
\cite{grothendieck}.

Let $\phi\in \Hom_k(A,I)$. Let $C^\phi\in \Hom_k(A\otimes_k A,I)$ be
defined by
\[ C^\phi(x,y)=x\phi(y)-\phi(xy)+\phi(x)y.\]
One checks that $C^\phi\in exan_k(A,I)$ for all $\phi\in \Hom_k(A,I)$.

\begin{definition} Let $\exan_k^{inn}(A,I)$ be the subset of
  $\exan_k(A,I)$ of maps $C^\phi$ for $\phi\in \Hom_k(A,I)$.
\end{definition}

\begin{lemma} The set $\exan_k^{inn}(A,I)\subseteq \exan_k(A,I)$ is a
  left and right sub $C(A)$-module.
\end{lemma}
\begin{proof} The proof is left to the reader as an exercise.
\end{proof}

\begin{definition} Let $\Exan_k^{inn}(A,I)\subseteq \Exan_k(A,I)$ be
  the image of $\exan_k^{inn}(A,I)$ under the bijection
  $\exan_k(A,I)\cong \Exan_k(A,I)$.
\end{definition}

It follows $\Exan_k^{inn }(A,I)\subseteq \Exan_k(A,I)$ is a left and
right sub $C(A)$-module.

Recall the definition of the \emph{Hochschild complex}:
\begin{definition} Let $A$ be an associative $k$-algebra and let $I$
  be a left and right $A$-module. Let $C^p(A,I)=\Hom_k(A^{\otimes
    p},I)$.  Let $d^p:C^p(A,I)\rightarrow C^{p+1}(A,I)$ be defined as
  follows:
\[ d^p(\phi)(a_1\otimes \cdots \otimes a_{p+1})=a_1\phi(a_2\otimes
\cdots \otimes a_{p+1})+\]
\[ \sum_{1\leq i \leq p}(-1)^i\phi(a_1\otimes \cdots \otimes
a_ia_{i+1}\otimes \cdots \otimes a_{p+1})+(-1)^{p+1}\phi(a_1\otimes
\cdots \otimes a_p)a_{p+1}.\]
We let $\HH^i(A,I)$ denote the i'th cohomology of this complex.
It is the \emph{$i$'th Hochschild cohomology} of $A$ with values in
$I$.
\end{definition}

\begin{proposition} \label{exact} There is an exact sequence
\[ 0\rightarrow \Exan_k^{inn}(A,I) \rightarrow \Exan_k(A,I)\rightarrow
\HH^2(A,I)\rightarrow 0 \]
of left and right $C(A)$-modules.
\end{proposition}
\begin{proof} The proof is left to the reader as an exercise.
\end{proof}

\begin{example} Characteristic classes of $L$-connections.
\end{example}
Let $A$ be a commutative $k$-algebra and let $\alpha:L\rightarrow
\Der_k(A)$ be a Lie-Rinehart algebra. Let $W$ be a left $A$-module
with an $L$-connection $\nabla:L\rightarrow \End_k(W)$. In
\cite{maa12} we define a characteristic class $c_1(E)\in
\H^2(L|_U,\O_U)$ when $W$ is of finite presentation, $U\subseteq
\spec(A)$ is the open set where $W$ is locally free and
$\H^2(L|_U,\O_U)$ is the Lie-Rinehart cohomology of $L|_U$ with values
in $\O_U$. If $L$ is locally free it follows $\H^2(L,A)\cong
\Ext_{U(L)}^2(A,A)$ where $U(L)$ is the generalized universal
enveloping
algebra of $L$. There is an obvious structure of left and right
$U(L)$-module on $\End_k(A)$ and an isomorphism
\[ \HH^2(U(L),\End_k(A))\cong \Ext_{U(L)}^2(A,A)\]
of abelian groups. The exact sequence \ref{exact} gives a sequence
\[ 0\rightarrow \Exan_k^{inn}(U(L), \End_k(A)) \rightarrow
\Exan_k(U(L),\End_k(A))\rightarrow \]
\[  \Ext^2_{U(L)}(A,A) \rightarrow 0\]
with $A=U(L)$ and $I=\End_k(A)$. If we can construct a lifting
\[ \tilde{c}_1(W)\in \Exan_k(U(L),\End_k(A)) \]
of the class
\[ c_1(W)\in \Ext_{U(L)}^2(A,A)=\HH^2(U(L),\End_k(A)) \]
we get a generalization of the
characteristic class from \cite{maa12} to arbitrary Lie-Rinehart
algebras $L$. This problem will be studied in a future paper on the
subject (see \cite{maa4}).

\begin{example} Non-commutative Kodaira-Spencer maps \end{example}

Let $A$ be an associative $k$-algebra and let $M$ be a left
$A$-module. Let $D^1(A)\subseteq \End_k(A)$ be the \emph{module of
  first
order differential operators} on $A$. It is defined as follows:
An element $\pa\in \End_k(A)$ is in $D^1(A)$ if and only if
$[\pa,a]\in D^0(A)=A\subseteq \End_k(A)$ for all $a\in A$. Define the
following map:
\[ f:D^1(A)\rightarrow \Hom_k(A, \End_k(M)) \]
by
\[ f(\pa)(a,m)=[\pa,a]m=(\pa(a)-a\pa(1))m.\]
Here $\pa\in D^1(A), a\in A$ and $m\in M$. Since $[\pa,a]\in A$ we get
a well defined map. Let for any $a\in A$ and $m\in M$ $\phi_a(m)=am$. 
It follows $\phi_a\in \End_k(M)$ is an endomorphism of $M$.
We get
\[ f(\pa)(ab,m)=(\pa(ab)-ab\pa(1))m=(\pa
\phi_{ab}-\phi_{ab}\pa)(1)m=\]
\[ (\pa \phi_{ab}-\phi_a \pa \phi_b +\phi_a \pa \phi_b
-\phi_{ab}\pa)(1)m=\]
\[(\pa \phi_a -\phi_a\pa)\phi_b(1)m+\phi_a(\pa
\phi_b-\phi_b\pa)(1)m=\]
\[ f(\pa)(a,bm)+af(\pa)(b,m).\]
Hence
\[ f(\pa)(ab)=af(\pa)(b)+f(\pa)(a)b \]
for all $\pa\in D^1(A)$ and $a,b\in A$.
The Hochschild complex gives a map
\[ d^1:\Hom_k(A,\End_k(M))\rightarrow \Hom_k(A\otimes A,\End_k(M)) \]
and
\[ ker(d^1)=\Der_k(A,\End_k(M)).\]
It follows we get a map
\[ f:D^1(A)\rightarrow \Der_k(A, \End_k(M)).\]
We get an induced map
\[ f:D^1(A)\rightarrow \HH^1(A, \End_k(M))=\Ext^1_A(M,M) .\]

\begin{lemma}  The following holds: $f(D^0(A))=f(A)=0$
\end{lemma}
\begin{proof} The proof is left to the reader as an exercise.
\end{proof}

One checks that $D^1(A)/D^0(A)=D^1(A)/A\cong \Der_k(A)$.
It follows we get an induced map
\[ g:\Der_k(A)=D^1(A)/D^0(A) \rightarrow \Ext^1_A(M,M) \]
the \emph{non-commutative Kodaira-Spencer map}.

\begin{lemma} Assume $A$ is commutative. The following holds:
\begin{align}
&\label{lr1}\text{$\mathbb{V}_M=ker(g)\subseteq \Der_k(A)$ is a
  Lie-Rinehart algebra.}\\
&\label{lr2}\text{$g(\delta)=0 \iff \exists \phi\in \End_k(M),
  \phi(am)=a\phi(m)+\delta(a)m.$}\\
&\label{lr3}\text{$\exists \nabla\in \Hom_k(\mathbb{V}_M, \End_k(M))$
  with
  $\nabla(\delta)(am)=a\nabla(\delta)(m)+\delta(a)m.$}\\
&\label{lr4}\text{$\mathbb{V}_M$ is the maximal Lie-Rinehart algebra
satisfying \ref{lr3}.}
\end{align}
\end{lemma}
\begin{proof} We first prove \ref{lr1}: Assume
  $g(\delta)=g(\eta)=0$.  By definition this is if and only if there
  are maps $\phi,\psi\in \End_k(M)$ such that the following holds:
\begin{align}
&\label{a1} d^0\phi=g(\delta) \\
&\label{a2}d^0\psi=g(\eta).
\end{align}
One checks that condition \ref{a1} and \ref{a2} hold if and only if the
following hold:
\[ \phi(am)=a\phi(m)+\delta(a)m \]
and
\[ \psi(am)=a\psi(m)+\eta(a)m.\]
We claim : $d^0[\delta,\eta]=g([\delta,\eta])$:
We get
\[ [\phi,\psi](am)=\phi\psi(am)-\psi\phi(am)=\]
\[\phi(a\psi(m)+\eta(a)m)-\psi(a\phi(m)+\delta(a)m)=\]
\[a\phi\psi(m)+\delta(a)\psi(m)+\eta(a)\phi(m)+\delta\eta(a)m -\]
\[ a\psi\phi(m)-\eta(a)\phi(m)-\delta(a)\psi(m)-\eta\delta(a)m= \]
\[a[\phi,\psi](m)+[\delta,\eta](a)m.\]
Hence $g([\delta,\eta])=0$ and $\mathbb{V}_M\subseteq \Der_k(A)$ is a
$k$-Lie algebra. It is an $A$-module since $g$ is $A$-linear, hence it
is a Lie-Rinehart algebra. Claim \ref{lr1} is proved.
Claim \ref{lr2} and \ref{lr3} follows from the proof of
\ref{lr1}. Claim \ref{lr4} is obvious and the Lemma is proved.
\end{proof}

The Lie-Rinehart algebra $\mathbb{V}_M$ is the \emph{linear
  Lie-Rinehart algebra} of $M$.

Let in the following $E$ be a left and right $A$-module.

\begin{definition} Let
\[ \Pr^1_I(E)=I\otimes_A E\oplus E \]
be the \emph{first order $I$-jet bundle of $E$}.
\end{definition}

Pick a derivation $d\in \Der_k(A,I)$ of left and right modules. This
means
\[ d(xy)=xd(y)+d(x)y \]
for all $x,y\in A$.
Let $B^C=I\oplus^C A$ and define the following left $B^C$-action on
$\Pr^1_I(E)$:
\[ (u,x)(w\otimes e,f)=(u\otimes f+xw\otimes e+d(x)\otimes f, xf) \]
for any elements $(u,x)\in B^C$ and $(w\otimes e,f)\in \Pr^1_I(E)$.

\begin{proposition} The abelian group $\Pr^1_I(E)$ is a left
  $B^C$-module if and only if $C(y,x)\otimes f=0$ for all $y,x\in A$
  and $f\in E$.
\end{proposition}
\begin{proof} One easily checks that for any $a,b\in B^C$ and $l,j\in
  \Pr^1_I(E)$ the following hold:
\[ (a+b)i=ai+bi \]
\[ a(i+j)=ai+aj.\]
Moreover
\[ \one i=i.\]
It remains to check that $a(bi)=(ab)i$. Let $a=(v,y)\in B^C$ and
$b=(u,x)\in B^C$. Let also $i=(w\otimes e, f)\in \Pr^1_I(E)$.
We get
\[ a(bi)=(v,y)((u,x)(w\otimes e,f))=(vx\otimes f+yu\otimes
f+yxw\otimes e +d(yx)\otimes f, yxf).\]
We also get
\[ (ab)i=(vx\otimes f+yu\otimes f +yxw\otimes e+d(yx)\otimes
f+C(y,x)\otimes f,yxf).\]
It follows that
\[ (ab)i-a(bi)=0\]
if and only if
\[ C(y,x)\otimes f=0 ,\]
and the claim of the Proposition follows.
\end{proof}

Note the abelian group $\Pr^1_I(E)$ is always a left $A$-module and
there is an exact sequence of left $A$-modules
\[ 0\rightarrow I\otimes E \rightarrow \Pr^1_I(E) \rightarrow E
\rightarrow 0 \]
defining a characteristic class
\[ c_I(E)\in \Ext^1_A(E,E\otimes I).\]
The class $c_I(E)$ has the property that $c_I(E)=0$ if and only if $E$
has an $I$-connection:
\[ \nabla:E\rightarrow I\otimes E \]
with
\[ \nabla(xe)=x\nabla(e)+d(x)\otimes e.\]
Let $J\subseteq I\subseteq B^C$ be the smallest two sided ideal
containing $Im(C)$
where $C:A\otimes_k A\rightarrow I$ is the cocycle defining $B^C$.
Let $D^C=B^C/J$ and $I^C=I/J$. We get a square zero extension
\[ 0\rightarrow I^C \rightarrow D^C \rightarrow A\rightarrow 0 \]
of $A$ by the square zero ideal $I^C$. It follows $D^C=I^C\oplus A$ as
abelian group. Since $\overline{C(x,y)}=0$ in $I^C$ it follows $D^C$
has
a well defined  associative multiplication defined by
\[ (u,x)(v,y)=(uy+xv, xy).\]
Also $D^C$ is the largest quotient of $B^C$ such that the ring
homomorphism $B^C\rightarrow D^C$ fits into a commutative diagram of
square zero extensions
\[
\diagram 0 \rto & I \rto \dto & B^C \rto \dto & A \rto \dto^= & 0 \\
         0 \rto & I^C \rto  & D^C \rto  & A \rto  & 0 .
\enddiagram
\]

\begin{definition} Let
\[ \Pr^1_{I^C}(E)=I^C\otimes E\oplus E \]
be the \emph{first order $I^C$-jet bundle of} $E$.
\end{definition}

\begin{example} First order commutative jets.
\end{example}
Let $k\rightarrow A$ be a commutative $k$-algebra and let $I\subseteq
A\otimes_k A$ be the ideal of the diagonal. Let $\Pr_A^1=A\otimes
A/I^2$ and $\Omega_A^1=I/I^2$. We get an exact sequence of left
$A$-modules
\begin{align}
\label{seq}0 \rightarrow \Omega^1_A \rightarrow \Pr^1_A \rightarrow A
\rightarrow 0 .
\end{align}
It follows $\Pr^1_A\cong \Omega^1_A \oplus A$ with the following
product:
\[ (\omega,a)(\eta,b)=(\omega a+b\eta,ab)\]
hence the sequence \ref{seq} splits.
Let $\Pr_A^1(E)=\Omega^1_A\otimes E\oplus E$ be the first order
$\Omega^1_A$-jet of $E$. We get an exact sequence of left $A$-modules
\[ 0\rightarrow \Omega_A^1\otimes E\rightarrow \Pr^1_A(E)\rightarrow
E\rightarrow 0 .\]
Since the sequence \ref{seq} splits it follows $\Pr_A^1(E)$ is a
lifting of $E$ to the first order jet $\Pr_A^1$.

\section{Atiyah classes and Kodaira-Spencer classes}

In this section we define and prove some properties of Atiyah classes
and Kodaira-Spencer classes.  

Let $X$ be any scheme defined over an arbitrary basefield $F$ and let
$\Pic(X)$ be the \emph{Picard group} of $X$. Let $\O^*\subseteq \O_X$
be the following subsheaf of abelian groups: For any open set
$U\subseteq X$ the group $\O(U)^*$ is the multiplicative group of
units in $\O_X(U)$. Define for any open set $U\subseteq X$ the
following morphism:
\[ \dlog: \O(U)^* \rightarrow \Omega^1_X(U) \]
defined by
\[ \dlog(x)=d(x)/x ,\]
where $d$ is the universal derivation and $x\in \O(U)^*$. 

\begin{lemma} The following hold:
\[ \dlog(xy)=\dlog(x)+\dlog(y)\]
for $x,y\in \O(U)^*$
\end{lemma}
\begin{proof} The proof is left to the reader as an exercise.
\end{proof}

Hence $\dlog:\O^*\rightarrow \Omega^1_X$ defines a map of sheaves of abelian
groups. The map $\dlog$ induce a map on cohomology
\[ \dlog: \Pic(X)=\H^1(X,\O^*)\rightarrow \H^1(X,\Omega^1_X) \]
and by definition
\[ c_l(\L)=\dlog(\L) \in \H^1(X, \Omega^1_X). \]

Let $\I\subseteq \Omega^1_X$ be any sub $\O_X$-module and let
$\F=\Omega^1_X/\I$ be the quotient sheaf. We get a derivation
\[ d:\O_X\rightarrow \F \]
by composing with the universal derivation. We get a canonical map
\[ \H^1(X,\Omega^1_X)\rightarrow \H^1(X,\F) \]
and we let
\[ \overline{c}_1(\L)\in \H^1(X,\F) \]
be the image of $c_1(\L)$ under this map.

\begin{definition} The class $c_1(\L)\in \H^1(X,\Omega^1_X)$ is the 
\emph{first Chern class} of the line bundle $\L\in \Pic(X)$. The
class $\overline{c}_1(\L)\in \H^1(X,\F)$ is the \emph{generalized
  first Chern class} of $\L$.
\end{definition}

Let $\E$ be any
$\O_X$-module
and consider the following sequence of sheaves of abelian groups:
\[ 0\rightarrow \F\otimes \E \rightarrow \Pr^1_{\F}(\E)\rightarrow \E
\rightarrow 0 \]
where
\[ \Pr^1_{\F}(\E)=\F\otimes \E\oplus \E \]
as sheaf of abelian groups. Let $s$ be a local section of $\O_X$ and
let
$(x\otimes e, f)$ be a local section of $\Pr^1_{\F}(\E)$ over some
open set $U$. Make the following definition:
\[ s(x\otimes e,f)=(sx\otimes e+ds\otimes f, sf).\]
It follows the sequence
\[ 0\rightarrow \F\otimes \E \rightarrow \Pr^1_{\F}(\E)\rightarrow \E
\rightarrow 0 \]
is a short exact sequence of sheaves of abelian groups.
It is called the \emph{Atiyah-Karoubi sequence}.
\begin{definition} An $\F$-connection $\nabla$ is a map
\[\nabla:\E\rightarrow \F\otimes \E \]
of sheaves of abelian groups with 
\[ \nabla (se)=s\nabla(e)+d(s)\otimes e.\]
\end{definition}

\begin{proposition} The Atiyah-Karoubi sequence is an exact sequence of
  left $\O_X$-modules. It is left split by an $\F$-connection.
\end{proposition}
\begin{proof} We first show it is an exact sequence of left $\O_X$-modules. The $\O_X$-module
  structure is twisted by the derivation $d$, hence we must verify
  that this gives a well defined left $\O_X$-structure on $\Pr^1_{\F}(\E)$.
Let $\omega=(x\otimes e,f)$ be a local section of $\Pr^1_{\F}(\E)$ and
let $s,t$ be local sections of $\O_X$. We get the following
calculation:
\[ (st)\omega=(st)(x\otimes e,f)=((st)x\otimes e+d(st)\otimes f, (st)f)=\]
\[(stx\otimes e+sdt\otimes f+(ds)t\otimes f, stf)=(s(tx\otimes
e+dt\otimes f)+ds\otimes tf, s(tf)) =\]
\[s(tx\otimes e+dt\otimes f,tf)=s(t(x\otimes e,f))=s(t\omega) .\]
It follows $\Pr^1_{\F}(\E)$ is a left $\O_X$-module and the sequence
is left exact. Assume
\[s:\E \rightarrow \Pr_{\F}(\E)=\F\otimes \E\oplus \E\]
is a left splitting. It follows $s(e)=(\nabla(e),e)$ for $e$ a local
section of $\E$. It follows $\nabla$ is a generalized connection and
the Theorem is proved.
\end{proof}

Note: If $\I=0$ we get $\Pr^1_{\F}(\E)=\Pr^1_X(\E)$ is the first order
jet bundle of $\E$ and the exact sequence above specialize to the well
known \emph{Atiyah sequence}:
\[ 0\rightarrow \Omega^1_X\otimes \E\rightarrow \Pr^1_X(\E)
\rightarrow \E \rightarrow 0.\]
The Atiyah sequence is left split by a connection
\[ \nabla: \E\rightarrow \Omega^1_X\otimes \E.\]
The $\O_X$-module $\Pr^1_{\F}(\E)$ is the \emph{generalized first
  order jet bundle of $\E$}.

\begin{definition} The characteristic class
\[ \AT(\E)\in \Ext^1_{\O_X}(\E, \F\otimes \E) \]
is called the \emph{Atiyah class} of $\E$.
\end{definition}
The class $\AT(\E)$ is defined for an arbitrary $\O_X$-module $\E$ and an arbitrary sub module
$\I\subseteq \Omega^1_X$.

Assume $\E=\L\in \Pic(X)$ is a line bundle on $X$. We get isomorphisms
\[ \Ext_{\O_X}(\L, \L\otimes \F)\cong \Ext^1_{\O_X}(\O_X,
\L^*\otimes \L\otimes \F)\cong \]
\[ \Ext^1_{\O_X}(\O_X, \F)\rightarrow  \H^1(X,\F).\]
We get a morphism
\[ \phi: \Ext^1_{\O_X}(\L,\L\otimes \F)\rightarrow \H^1(X,\F).\] 

\begin{proposition} The following hold:
\[\phi(\AT(\L))=\overline{c}_1(\L).\]
Hence the Atiyah class calculates the generalized first Chern class of a line bundle.
\end{proposition}
\begin{proof} Let $\I=0$. It is well known that $\AT(\L)$ calculates
  the first Chern class $c_1(\L)$. From this the claim of the
  Proposition follows.
\end{proof}

Let $T_X$ be the tangent sheaf of $X$. It has the property that for
any open affine set $U=\spec(A)\subseteq X$ the local sections
$T_X(U)$ equals the module $\Der_F(A)$ of derivations of $A$.
Let $\V\subseteq T_X$ be the subsheaf of local sections 
$\partial$ of $T_X$ with the following property: The section
$\partial \in T_X(U)$ lifts to a local
section $\nabla(\partial)$ of $\End_F(\E|_U)$ with the following property:
\[ \nabla(\partial):\E|_U\rightarrow \E|_U \]
satisfies
\[ \nabla(\partial)(se)=s\nabla(\partial)(e)+\partial(s)e.\]
It follows $\V\subseteq T_X$ is a sheaf of Lie-Rinehart algebras - the
\emph{Kodaira-Spencer sheaf} of $\E$. 

Define for any local sections $a,b$ of $\O_X$,$\partial $ of $\V$ and
$e$ of $\E$ the following:
\[ L(a,\partial)(e)=a\nabla(\partial)(e)-\nabla(a\partial)(e).\]

\begin{lemma} It follows $L(a,\partial)\in \End_{\O_U}(\E|_U)$ .
\end{lemma}
\begin{proof} The following hold:
\[ L(a,\partial)(be)=a\nabla(\partial)(be)-\nabla(a\partial)(be)=\]
\[
a(b\nabla(\partial)(e)+\partial(b)e)-b\nabla(a\partial)(e)-a\partial(b)e
=\]
\[ab\nabla(\partial)(e)+a\partial(b)e-b\nabla(a\partial)(e)-a\partial(b)e=\]
\[b(a\nabla(\partial)(e)-\nabla(a\partial)(e))=b(a\nabla(\partial)-\nabla(a\partial))(e)=\]
\[ bL(a,\partial)(e)\]
and the Lemma is proved.
\end{proof}

\begin{lemma} The following formula hold:
\[ L(ab,\partial)=aL(b,\partial)+L(a,b\partial) \]
for all local sections $a,b$ and $\partial$.
\end{lemma}
\begin{proof} We get
\[ L(ab,\partial)=ab\nabla(\partial)-\nabla(ab\partial)=\]
\[
ab\nabla(\partial)-a\nabla(b\partial)+a\nabla(b\partial)-\nabla(ab\partial)=\]
\[ a(b\nabla(\partial)-\nabla(b\partial))+(a\nabla-\nabla
a)(b\partial)=\]
\[aL(b,\partial)+L(a,b\partial),\]
and the Lemma is proved.
\end{proof}

Let $\LR(\V)=\End_{\O_X}(\E)\oplus \V$ be the \emph{linear
  Lie-Rinehart algebra} of $\E$. Let $\LR(\V)$ have the following left
$\O_X$-module structure:
\[ a(\phi, \partial)=(a\phi+L(a,\partial), a\partial).\]
Here $a,\phi$ and $\partial$ are local sections of $\O_X$,
$\End_{\O_X}(\E)$ and $\V$. We twist the trivial $\O_X$ structure on
$\End_{\O_X}(\E)\oplus \V$ with the element $L$.
We get a sequence of sheaves of abelian groups
\[ 0\rightarrow \End_{\O_X}(\E)\rightarrow^i \LR(\V) \rightarrow^p \V \rightarrow 0\]
where $i$ and $p$ are the canonical maps.
An $\O_X$-linear  map
\[ \nabla:\V \rightarrow \End_F(\E) \]
satisfying
\[ \nabla(\partial)(ae)=a\nabla(\partial)(e)+\partial(a)e \]
is a \emph{$\V$-connection} on $\E$.

\begin{proposition}  The sequence defined above is an exact sequence of
  left $\O_X$-modules. It is left split by a $\V$-connection $\nabla$.
\end{proposition}
\begin{proof} We need to check that $\LR(\V)$ has a well defined left
  $\O_X$-module structure. By definition
\[ a(\phi, \partial)=(a\phi+L(a,\partial),a\partial).\]
We get
\[(ab)x= (ab)(\phi,\partial)=((ab)\phi+L(ab,\partial), (ab)\partial)=\]
\[(ab\phi+aL(b,\partial)+L(a,b\partial),ab\partial)=\]
\[a(b\phi+L(b,\partial), b\partial)=a(b(\phi,\partial))=a(bx)\]
and it follows the sequence is  a left exact sequence of
$\O_X$-modules. If
\[ s:\V \rightarrow \End_{\O_X}(\E)\oplus \V =\LR(\V) \]
is a section it follows $s(e)=(\nabla(e),e)$. One checks that
$\nabla$ is a $\V$-connection, and the Theorem is proved.
\end{proof}

\begin{definition} We get a characteristic class
\[ \KS(\E)\in \Ext^1_{\O_X}(\V, \End_{\O_X}(\E)) \]
the \emph{Kodaira-Spencer class} of $\E$.
\end{definition}

Assume $\V$ is locally free and $\E=\L\in \Pic(X)$ is a line bundle on
$X$. Assume also $\V^*=\F=\Omega^1_X/\I$ for some submodule $\I$. 
We get the following calculation:

\[ \Ext^1_{\O_X}(\V,\End_{\O_X}(\L))\cong \Ext^1_{\O_X}(\O_X,
\End_{\O_X}(\L)\otimes \V^*)\cong \]
\[ \Ext^1_{\O_X}(\O_X, \End_{\O_X}(\L) \otimes \F)\rightarrow
\H^1(X,\F).\]
We get a map
\[ \psi:\Ext^1_{\O_X}(\V,\End_{\O_X}(\L))\rightarrow \H^1(X,\F)\]
of sheaves.

\begin{proposition} The following hold: There is an equality
\[ \psi(\KS(\L))=\overline{c}_1(\L) \]
in $\H^1(X,\F)$. 
Hence the Kodaira-Spencer class calculates the class $\overline{c}_1(\L)$.
\end{proposition}
\begin{proof}The proof is left to the reader as an exercise.
\end{proof}

We get the following diagram expressing the relationship between the
characteristic classes defined above:
\[
\diagram \Ext^1_{\O_X}(\VL, \End_{\O_X}(\L)) \drto^{\psi} & & \\
                                &  \H^1(X,\F) &      \Pic(X) \lto_{\overline{c}_1(-)} \\
     \Ext^1_{\O_X}(\L, \F\otimes \L) \urto^{\phi} & & 
\enddiagram 
\] 
The following equation holds in $\H^1(X,\F)$:
\[ \phi(\AT(\L))=\psi(\KS(\L))=\overline{c}_1(\L) .\]


\begin{thebibliography}{4}

\bibitem{andre} M. Andre, Homologie des algebres commutatives,
  \emph{Grundlehren Math. Wiss.} no. 206 (1974)

\bibitem{atiyah} M. Atiyah, Complex analytic connections in fibre
  bundles, \emph{Trans. AMS} no. 85 (1957)

\bibitem{grothendieck} A. Grothendieck, EGA IV Etude locale de schemas
  et des morphismes de schemas, \emph{Publ. Math. IHES} no. 20 (1964)

\bibitem{karoubi} M. Karoubi, Homologie cyclique et K-theorie,
  \emph{Asterisque} no 149 (1987)

\bibitem{maa1} H. Maakestad, A note on the principal parts on
  projective space and linear representations,
\emph{Proc. of the AMS} Vol. 133 no. 2 (2004) 

\bibitem{maa12} H. Maakestad, Chern classes and Lie-Rinehart algebras,
  \emph{Indagationes Math.} (2008)

\bibitem{maa4} H. Maakestad, Chern classes and $\Exan$ functors, \emph{In progress} (2009)

\bibitem{maa2} H. Maakestad, Principal parts on the projective line
  over arbitrary rings, \emph{Manuscripta Math.} 126, no. 4 (2008) 




\end{thebibliography}
\end{document}